\documentclass[12pt,fleqn]{article}
\usepackage{mathrsfs}
\usepackage{amsfonts}
\usepackage{amssymb}
\usepackage{amsthm}
\usepackage{enumerate}
\usepackage{latexsym,amsmath}
\usepackage{graphicx}
\usepackage{abstract}
\usepackage{amsmath}
\usepackage{caption,subcaption}
\title{On the Laplacian spectra of some double join operations of graphs}

\date{}

\begin{document}

\author{Gui-Xian Tian$^a$\footnote{Corresponding author. E-mail: gxtian@zjnu.cn or guixiantian@163.com(G.-X. Tian)}, Jing-Xiang He$^a$, Shu-Yu Cui$^b$,\\
{\small{\it $^a$College of Mathematics, Physics and Information Engineering,}}\\
{\small{\it Zhejiang Normal University, Jinhua, Zhejiang, 321004,
P.R. China}} \\
{\small{\it $^b$Xingzhi College, Zhejiang Normal University, Jinhua,
Zhejiang, 321004, P.R. China}} }\maketitle

\begin{abstract}
Many variants of join operations of graphs have been introduced and
their spectral properties have been studied extensively by many
researchers. This paper mainly focuses on the Laplacian spectra of
some double join operations of graphs. We first introduce the
conception of double join matrix and provide a complete information
about its eigenvalues and the corresponding eigenvectors. Further,
we define four variants of double join operations based on
subdivision graph, $Q$-graph, $R$-graph and total graph. Applying
the result obtained for the double join matrix, we give an explicit
complete characterization of the Laplacian eigenvalues and the
corresponding eigenvectors of four variants in terms of the
Laplacian eigenvalues and the eigenvectors of the factor graphs.
These results generalize some well-known results about some join
operations of graphs.



\end{abstract}

\section*{1. Introduction}

\indent\indent Throughout this paper, all graphs considered are
finite simple graphs. Let $G=(V,E)$ be a graph with vertex set
$V=\{v_{1},v_{2},\ldots,v_{n}\}$ and edge set $E(G)$. The adjacency
matrix $A(G)=(a_{ij})$ of $G$ is an $n\times n$ matrix where
$a_{ij}=1$ whenever $v_{i}$ and $v_{j}$ are adjacent in $G$ and
$a_{ij}=0$ otherwise. The degree of $v_{i}$ in $G$ is denoted by
$d_{i}=d_{G}(v_{i})$. Let $D(G)$ be the degree diagonal matrix of
$G$ with diagonal entries $d_{1},d_{2},\ldots,d_{n}$. The Laplacian
matrix $L(G)$ of $G$ is defined as $D(G)-A(G)$. The signless
Laplacian matrix of $G$ is defined as $|L|(G)=D(G)+A(G)$. For an
$n\times n$ matrix $M$ associated to $G$, the set of all the
eigenvalues of $M$ is called the spectrum of matrix $M$ or graph
$G$. In particular, if $M$ is the adjacency matrix $A(G)$ of $G$,
then the adjacency spectrum of $G$ is denoted by $\sigma(A(G)) =
(\nu_1 (G),\nu_2 (G), \ldots ,\nu_n (G)), $ where $\nu_1 (G)\leq
\nu_2 (G)\leq\cdots\leq \nu_n (G)$ are the eigenvalues of $A(G)$. If
$M$ is the Laplacian matrix $L(G)$ of $G$, then the Laplacian
spectrum of $G$ is denoted by $ \sigma(L(G)) = (\lambda_1
(G),\lambda_2 (G), \ldots ,\lambda_n (G)), $ where $\lambda_1
(G)\leq \lambda_2 (G)\leq\cdots\leq \lambda_n (G)$ are the
eigenvalues of $L(G)$. If $M$ is the signless Laplacian matrix
$|L|(G)$ of $G$, then the signless Laplacian spectrum of $G$ is
denoted by $ \sigma(|L|(G)) = (q_1 (G),q_2 (G), \ldots ,q_n (G)), $
where $q_1 (G)\leq q_2 (G)\leq\cdots\leq q_n (G)$ are the
eigenvalues of $|L|(G)$. For more review about the adjacency
spectrum, Laplacian spectrum and signless Laplacian spectrum of $G$,
readers may refer to
\cite{Brouwer2011,Chung1997,Cvetkovic1980,Cvetkovic2007,Dam2003,Grone1990,Merris1994}
and the references therein.

Determining the spectra of many graph operations is a basic and very
meaningful work in spectral graph theory. Up till now, many graph
operations such as Cartesian product, Kronecker product, graph with
$k$ (edge)-pockets, corona, edge corona, some variants of
(edge)corona, join, some variants of join have been introduced and
the adjacency spectra (Laplacian spectra, signless Laplacian spectra
as well) of these graph operations have also been determined in
terms of the corresponding spectra of the factor graphs in
\cite{Barik2008,Barik2015,Barik2007,Barik2017,Cardoso2013,Cui2012,Estrada2017,Hou2010,Indulal2012,Liu2017,Nath2014,Zhang2009}.
Moreover, it is known that the corresponding spectra of these graph
operations can be used to construct infinitely many pairs of
cospectral
graphs\cite{Barik2008,Barik2007,Barik2017,Cui2012,Dam2003,Liu2017,Nath2014},
infinitely families of integral graphs\cite{Barik2015,Indulal2012}
and to investigate many other properties of graphs, such as the
Kirchhoff index\cite{Liu2017,Liu2015,Zhang2009}, the number of
spanning trees\cite{Barik2017,Hou2010,Liu2017} and so on. This paper
focuses on the Laplacian spectra of four new variants of double join
operations based on subdivision graph, $Q$-graph, $R$-graph and
total graph. The following definitions come from
\cite{Cvetkovic1980}, which will be required to define our new
operations.

Let $G$ be a connected graph with $n$ vertices and $m$ edges. The
\emph{subdivision graph} $S(G)$ is the graph obtained by inserting a
new vertex into every edge of $G$. The \emph{$Q$-graph} $Q(G)$ is
the graph obtained by inserting a new vertex into every edge of $G$
and by adding edges between those inserted vertex which lie on
adjacent edges of $G$.  The \emph{$R$-graph} $R(G)$ is the graph by
adding a new vertex corresponding to each edge of $G$ and by adding
edges between each added vertex and the corresponding edge's
endpoints. The \emph{total graph} $T(G)$ is the graph whose vertex
set is the union of vertex set and edge set of $G$, and two vertex
of $T(G)$ is adjacent whenever two corresponding elements are
incident or adjacent; see Figure 1 for example.

\begin{figure}[h]
\centering
\includegraphics[width=14cm]{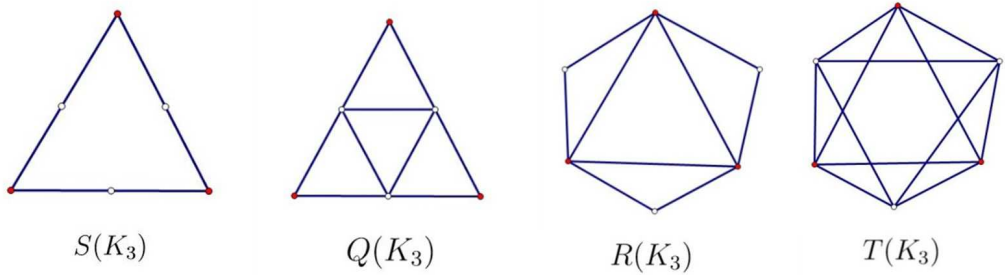}
\caption{$S(K_3)$, $Q(K_3)$, $R(K_3)$ and $T(K_3)$ for the complete
graph $K_3$. (Here new added vertices are white, the old vertices
are black.)}
\end{figure}

\paragraph{Definition 1.} Let $G$ be a connected graph with $n$ vertices and $m$
edges. Also let $G_1$ and $G_2$ be two graphs with $n_1$ and $n_2$
vertices, respectively. The \emph{subdivision double join}
$G^S\vee{(G^\bullet_1,G^\circ_2)}$ of $G$ , $G_1$ and $G_2$ is the
graph obtained from $S(G)$, $G_1$ and $G_2$ by joining every vertex
of $G$ to every vertex of $G_1$ and every vertex of $I(G)$ to every
vertex of $G_2$, where $I(G)$ denotes the vertex set of the added
new vertices in $S(G)$. Replaced $S(G)$ by $Q(G)$ ($R(G),T(G)$) in
this definition, the resulting graph is referred to as
\emph{$Q$-graph($R$-graph, total, respectively) double join} of
these graphs. Similarly, we denote them by $G^Q\vee{(G^\bullet_1 ,
G^\circ_2)}$, $G^R\vee{(G^\bullet_1 , G^\circ_2)}$ and
$G^T\vee{(G^\bullet_1 , G^\circ_2)}$, respectively.

\paragraph{Example 2.} Let $G$, $G_1$ and $G_2$ be the complete graph $K_3$, the paths $P_2$ and $P_3$, respectively. Figure
2 displays four graphs $K_3^S\vee{(P_2^\bullet, P_3^\circ)}$,
$K_3^Q\vee{(P_2^\bullet, P_3^\circ)}$, $K_3^R\vee{(P_2^\bullet,
P_3^\circ)}$ and $K_3^T\vee{(P_2^\bullet , P_3^\circ)}$ below.

\begin{figure}[h]
\centering
\includegraphics[width=15cm]{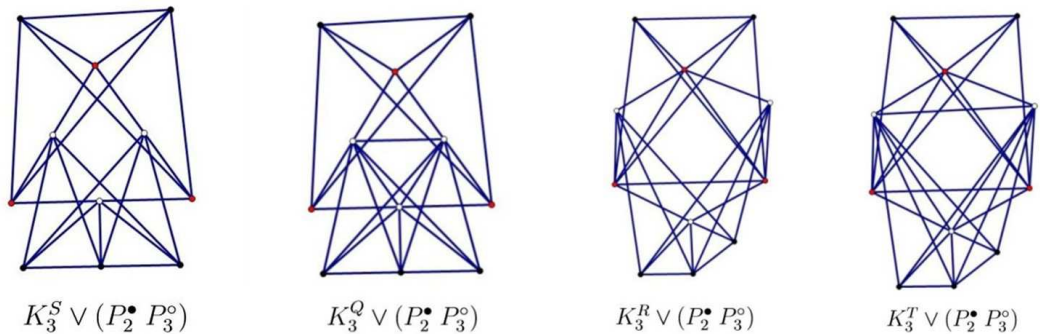}
\caption{Illustration of $K_3^S\vee{(P_2^\bullet, P_3^\circ)}$,
$K_3^Q\vee{(P_2^\bullet, P_3^\circ)}$, $K_3^R\vee{(P_2^\bullet,
P_3^\circ)}$ and $K_3^T\vee{(P_2^\bullet , P_3^\circ)}$ in Example
2.}
\end{figure}

Recently, many variants of join operations of graphs have been
introduced and their spectral properties have been studied by many
researchers. Cardoso et al.\cite{Cardoso2013} characterized
adjacency and Laplacian spectra of the $H$-join operation of graphs.
Estrada and Benzi\cite{Estrada2017} discussed the clustering,
assortativity and spectral properties of core-satellite graphs.
Remark that the core-satellite graph named in \cite{Estrada2017} is
a special join of some complete graphs. In \cite{Indulal2012}, the
adjacency spectra of the subdivision vertex(edge) joins of graphs
were computed in terms of the corresponding spectra of two regular
graphs. The author also constructed infinite family of new integral
graphs. Liu and Zhang\cite{Liu2017} determined the spectra,
(signless) Laplacian spectra of the subdivision vertex(edge) joins
for a regular graph $G_1$ and arbitrary graph $G_2$. As
applications, they constructed infinitely many pairs of cospectral
graphs and obtained the number of spanning trees and the Kirchhoff
index of the subdivision vertex(edge) joins. Remark that they
determined these spectra with the help of the coronal technique. But
this technique cannot describe completely the eigenvectors
corresponding to all the eigenvalues.

Motivated by these researches, we discuss the Laplacian spectra of
four new variants of double join operations based on subdivision
graph, $Q$-graph, $R$-graph and total graph namely,
$G^S\vee{(G^\bullet_1,G^\circ_2)}$, $G^Q\vee{(G^\bullet_1 ,
G^\circ_2)}$, $G^R\vee{(G^\bullet_1 , G^\circ_2)}$ and
$G^T\vee{(G^\bullet_1 , G^\circ_2)}$. The rest of this paper is
organized as follows. In Section 2, we shall introduce the
conception of double join matrix and provide a complete information
about its eigenvalues and the corresponding eigenvectors. In Section
3, applying the result obtained for the double join matrix, we give
an explicit complete characterization of the Laplacian spectra of
the four variants $G^S\vee{(G^\bullet_1,G^\circ_2)}$,
$G^Q\vee{(G^\bullet_1 , G^\circ_2)}$, $G^R\vee{(G^\bullet_1 ,
G^\circ_2)}$ and $G^T\vee{(G^\bullet_1 , G^\circ_2)}$ in terms of
the Laplacian spectra of the factor graphs. These results not only
generalize some well-known results, but also describe completely the
eigenvectors corresponding to all the Laplacian eigenvalues of these
graphs. In Section 4, we summarize our work and give some further
remarks.

\section*{2. Spectra of double join matrices}

\indent\indent Suppose that $A$, $C$, $D$, $E$ are real matrices of
order $p$, $q$, $r$, $s$, respectively and $B$ is a $p\times q$
matrix with $p<q$. Consider the following block matrix:
$$
\mathscr{D}_j=
   \begin{pmatrix}
   A & B& cJ_{p\times r}& \textbf{0}_{p\times s}\\
   B^T& C & \textbf{0}_{q\times r}& cJ_{q\times s}\\
   cJ_{r\times p}& \textbf{0}_{r\times q} & D & \textbf{0}_{r\times s}\\
   \textbf{0}_{s\times p}& cJ_{s\times q}& \textbf{0}_{s\times r}& E\\
   \end{pmatrix}  \notag,
$$
where $J_{p \times r} $ denotes the $p\times r$ matrix with every
entry is equal to 1 and $c=\pm 1$. Obviously, $\mathscr{D}_j$ is a
matrix of order $p+q+r+s$. Throughout  $\textbf{1}_n$ denotes the
column vector of order $n$ with every entry 1. We call
$\mathscr{D}_j$ the \emph{double join} matrix if $\mathscr{D}_j$
satisfies the following four conditions:
\begin{enumerate}[(i)]
    \item  If $\textbf{X}_i$ and $\textbf{Y}_i$ are the singular vector pairs of $B$
    corresponding to the singular values $b_i$  for $i=1,2,...,p$,
    then $\textbf{X}_i$ and $\textbf{Y}_i$ are the orthogonal unit eigenvectors of $A$ and $C$.
    Equivalently, $B \textbf{Y}_i=b_i \textbf{X}_i$ and $B^T \textbf{X}_i=b_i \textbf{Y}_i$ for $i=1,2,...,p$
    imply $A \textbf{X}_i=a_i \textbf{X}_i$, $C \textbf{Y}_i=c_i \textbf{Y}_i$
    where $a_i$ and $c_i$ are eigenvalues of matrices $A$ and $C$, respectively.
    \item Assume that $\textbf{X}_1=\frac{1}{\sqrt{p}}\textbf{1}_p$ and
    $\textbf{Y}_1=\frac{1}{\sqrt{q}}\textbf{1}_q$ are the unit
    eigenvectors of $A$ and $C$ corresponding to the eigenvalues $a_1 $ and $c_1$,
    respectively.
    \item If  $B\textbf{Y}_j=\textbf{0}_p$ for $j=p+1,p+2,...,q$,
    then $\textbf{Y}_j$ are the orthogonal eigenvectors of $C$, that is,
    $C \textbf{Y}_j=c_j\textbf{Y}_j$ for
    $j=p+1,p+2,...,q$, where $c_j$ are eigenvalues of $C$.
    \item Also assume that $\textbf{Z}_1=\frac{1}{\sqrt{r}}\textbf{1}_r$ and
    $\textbf{W}_1=\frac{1}{\sqrt{s}}\textbf{1}_s$ are the
    unit eigenvectors of $D$ and $E$ corresponding to the eigenvalues $d_1$ and $e_1$,
    respectively.
\end{enumerate}

Next, we shall give a full description of all the eigenvalues and
the corresponding eigenvectors for the double join matrix
$\mathscr{D}_j$.

\paragraph{Theorem 3.} The spectrum of the double join matrix
$\mathscr{D}_j$ consists of:
\begin{itemize}
  \item The eigenvalues $\lambda=d_i$ for $ i=2,3,...,r$;
  \item The eigenvalues $\lambda=e_i$ for $ i=2,3,...,s$;
  \item The eigenvalues $\lambda_{\pm}=\frac{a_i+c_i\pm \sqrt{(a_i-c_i)^2+4b_i^2}}{2}$ for $i=2,3,...,p$;
  \item The eigenvalues $\lambda=c_j$ for $j=p+1,p+2,...,q$;
  \item The remaining four eigenvalues  are given by the roots of the following
  equation:\\
  $\lambda^4-(e_1+c_1+a_1+d_1)\lambda^3+(e_1c_1-qs+(a_1+d_1)(e_1+c_1)+a_1d_1-pr-b^2_1)
  \lambda^2+((d_1+e_1)b^2_1-(a_1+d_1)(e_1c_1-qs)-(a_1d_1-pr)(e_1+c_1))\lambda+(a_1d_1-pr)(e_1c_1-qs)-d_1e_1b^2_1=0.$
\end{itemize}

\begin{proof}
Suppose that $\textbf{Z}_1,\textbf{Z}_2,\ldots,\textbf{Z}_r$ are the
orthogonal eigenvectors of $D$ corresponding to the eigenvalues
$d_1,d_2,...,d_r$, respectively. Firstly, consider the vectors
\begin{equation}
\;\;\;\;\;\;\;\;\;\;\;\;\;\;\;\;\;\;\;\;\;\;\;\;\;\;\;\;\;\;\;\;\;\textbf{x}=\begin{pmatrix}
\textbf{0}_p\\\textbf{0}_q\\\textbf{Z}_i\\\textbf{0}_s\\
\end{pmatrix}, \;\;\;\text{for}~ i=2,3,...,r.
\end{equation}
Notice that $J_{p\times r}\textbf{Z}_i=\textbf{0}_p$ as
$\textbf{Z}_i\bot \textbf{1}_r$. Then the equation
$\mathscr{D}_j\textbf{x}=\lambda \textbf{x}$ becomes
$$ \mathscr{D}_j\textbf{x}= \begin{bmatrix}
\textbf{0}_p\\\textbf{0}_q\\d_i\textbf{Z}_i\\\textbf{0}_s\\
\end{bmatrix}= \begin{bmatrix}
\textbf{0}_p\\\textbf{0}_q\\\lambda\textbf{Z}_i\\\textbf{0}_s\\
\end{bmatrix}.$$
So, $\lambda=d_i$ are the eigenvalues of double join matrix
$\mathscr{D}_j$ for $i=2,3...,r$.

Now suppose that $\textbf{W}_1,\textbf{W}_2,\ldots,\textbf{W}_s$ are
the orthogonal eigenvectors of $E$ corresponding to the eigenvalues
$e_1,e_2,...,e_s$, respectively. Next, consider the vectors
         \begin{equation}\label{2}
         \;\;\;\;\;\;\;\;\;\;\;\;\;\;\;\;\;\;\;\;\;\;\;\;\;\;\;\;\;\;\;\;\;\textbf{x}=\begin{pmatrix}
         \textbf{0}_p\\\textbf{0}_q\\\textbf{0}_r\\\textbf{W}_i\\
         \end{pmatrix},\;\;\; \text{for}~ i=2,3,...,s.
         \end{equation}
Then we plug (\ref{2}) into the equation
$\mathscr{D}_j\textbf{x}=\lambda \textbf{x}$. Notice that
$J_{q\times s}\textbf{W}_i=\textbf{0}_q$ as $\textbf{W}_i\bot
\textbf{1}_s$. Thus one gets
         $$ \mathscr{D}_j\textbf{x}= \begin{bmatrix}
         \textbf{0}_p\\\textbf{0}_q\\\textbf{0}_r\\e_i\textbf{W}_i\\
         \end{bmatrix}= \begin{bmatrix}
         \textbf{0}_p\\\textbf{0}_q\\\textbf{0}_r\\\lambda\textbf{W}_i\\
         \end{bmatrix}.$$
         So, $\lambda=e_i$ ($i=2,3,...,s$) are also eigenvalues of double join matrix $\mathscr{D}_j$.

Now consider the following vectors
         \begin{equation}\label{3}
         \;\;\;\;\;\;\;\;\;\;\;\;\;\;\;\;\;\;\;\;\;\;\;\;\;\;\;\;\;\;\;\;\;\textbf{x}=\begin{pmatrix}
         k_1\textbf{X}_i\\\textbf{Y}_i\\\textbf{0}_r\\\textbf{0}_s\\
         \end{pmatrix}, \;\;\; \text{for}~ i=2,3,...,p,
         \end{equation}
         where $k_1$ is an unknown constant to be determined. Notice that
$J_{r\times p}\textbf{X}_i=\textbf{0}_r$ and $J_{s\times
q}\textbf{Y}_i=\textbf{0}_s$ as $\textbf{X}_i\bot \textbf{1}_p$ and
$\textbf{Y}_i\bot \textbf{1}_q$. Again, plugging (\ref{3}) into the
equation $\mathscr{D}_j\textbf{x}=\lambda \textbf{x}$, we obtain
         $$\mathscr{D}_j\textbf{x}=\begin{bmatrix}
         k_1a_i\textbf{X}_i+b_i\textbf{X}_i\\k_1b_i\textbf{Y}_i+c_i\textbf{Y}_i\\\textbf{0}_r\\\textbf{0}_s\\
         \end{bmatrix}=\lambda \begin{bmatrix}
         k_1\textbf{X}_i\\\textbf{Y}_i\\\textbf{0}_r\\\textbf{0}_s\\
         \end{bmatrix},
         $$
         which reduces to the following conditions
         $$k_1a_i+b_i=\lambda k_1,\;\;\;k_1b_i+c_i=\lambda\notag.$$
         Eliminating $k_1$ from above conditions, one obtains
         \begin{gather}\label{4}
         \;\;\;\;\;\;\;\;\;\;\;\;\;\;\;\;\;\;\;\;\;\;\;\;\;\;\;\;\;\;\lambda^2-(a_i+c_i)\lambda+a_ic_i-b_i^2=0.
         \end{gather}
         The roots of the equation (\ref{4}) are
         $\lambda_{\pm}=\frac{(a_i+c_i)\pm\sqrt[2]{(a_i-c_i)^2+4b_i^2}}{2}$,
         which implies that the third part of theorem follows.

         Below, we consider the vectors
         \begin{equation}
         \;\;\;\;\;\;\;\;\;\;\;\;\;\;\;\;\;\;\;\;\;\;\;\;\;\;\;\;\;\;\;\;\;\textbf{x}=\begin{pmatrix}
         \textbf{0}_p\\\textbf{Y}_j\\\textbf{0}_r\\\textbf{0}_s\\
         \end{pmatrix},\;\;\; \text{for}~ j=p+1,p+2,...,q.
         \end{equation}
         Observe that $B\textbf{Y}_j=\textbf{0}_p$ and $J_{s\times
q}\textbf{Y}_j=\textbf{0}_s$ for $j=p+1,p+2,...,q$. Then the
equation $\mathscr{D}_j\textbf{x}=\lambda \textbf{x}$ becomes
         $$\mathscr{D}_j\textbf{x}=
         \begin{bmatrix}
         \textbf{0}_p\\c_j\textbf{Y}_j\\\textbf{0}_r\\\textbf{0}_s\\
         \end{bmatrix}=\lambda
         \begin{bmatrix}
         \textbf{0}_p\\\textbf{Y}_j\\\textbf{0}_r\\\textbf{0}_s\\
         \end{bmatrix}.
         $$
Hence, $\lambda =c_j$ ($j=p+1,p+2,...,q$) are eigenvalues of
$\mathscr{D}_j$. So far we have determined $p+q+r+s-4$ eigenvalues
of $\mathscr{D}_j$.

To determine the four remaining eigenvalues and the corresponding
eigenvectors, let
         \begin{equation}\label{6}
         \;\;\;\;\;\;\;\;\;\;\;\;\;\;\;\;\;\;\;\;\;\;\;\;\;\;\;\;\;\;\;\;\;\;\;\;\;\;\;\;\;\;\;\textbf{x}=\begin{pmatrix}
         k_1\textbf{1}_p\\k_2\textbf{1}_q\\k_3\textbf{1}_r\\\textbf{1}_s\\
         \end{pmatrix},
         \end{equation}
         where $k_1,k_2,k_3$ are three unknown constants to be
         determined. Note that
$B\textbf{1}_q=\sqrt{q}B\frac{\textbf{1}_q}{\sqrt{q}}=b_1\sqrt{\frac{q}{p}}\textbf{1}_p$
and $B^T\textbf{1}_p=b_1\sqrt{\frac{p}{q}}\textbf{1}_q$. Plugging
(\ref{6}) into the equation  $\mathscr{D}_j\textbf{x}=\lambda
\textbf{x}$, we get following conditions.

         $$ \begin{cases}
         k_1a_1+k_2b_1\sqrt{\frac{q}{p}}+ck_3r=\lambda k_1,\\
         k_1\sqrt{\frac{p}{q}}b_1+k_2c_1+cs=\lambda k_2,\\
         ck_1p+k_3d_1=\lambda k_3,\\
         ck_2q+e_1=\lambda.\\
         \end{cases}$$
         Eliminating $k_1,k_2$ and $k_3$ from above conditions, one
         obtains
         $$(\lambda-a_1-\frac{c^2pr}{\lambda-d_1})(\lambda-c_1-\frac{c^2qs}{\lambda-e_1})=b_1^2.$$
         Note $c^2=1$. We may reduce this equation to\\\\
         $\lambda^4-(e_1+c_1+a_1+d_1)\lambda^3+(e_1c_1-qs+(a_1+d_1)(e_1+c_1)+a_1d_1-pr-b^2_1)\lambda^2
         +((d_1+e_1)b^2_1-(a_1+d_1)(e_1c_1-qs)-(a_1d_1-pr)(e_1+c_1))\lambda+(a_1d_1-pr)(e_1c_1-qs)-d_1e_1b^2_1=0.$\\\\
The proof of this theorem is completed.
    \end{proof}

\section*{3. Laplacian spectra of double join operations of graphs}

\indent \indent In this section, applying the result obtained in
Theorem 3, we shall give an explicit complete characterization of
the Laplacian spectra of four variants of the double join operations
$G^S\vee{(G^\bullet_1,G^\circ_2)}$, $G^Q\vee{(G^\bullet_1 ,
G^\circ_2)}$, $G^R\vee{(G^\bullet_1 , G^\circ_2)}$ and
$G^T\vee{(G^\bullet_1 , G^\circ_2)}$ in terms of the Laplacian
spectra of the factor graphs.

We first focus on determining the Laplacian spectra of the
subdivision double join $G^S\vee{(G^\bullet_1, G^\circ_2)}$ for a
regular graph $G$ and two arbitrary graphs $G_1$, $G_2$.

\paragraph{Theorem 4.} Let $G$ be a $k$-regular graph with $n$ vertices and $m$
    edges. Also let $G_1$ and $G_2$ be two arbitrary graph with $n_1$ and $n_2$
    vertices, respectively. Then the Laplacian spectrum of  $G^S\vee{(G^\bullet_1
, G^\circ_2)}$ consists of:
    \begin{enumerate}[(i)]
        \item $\lambda_i(G_1)+n$, for $i=2,3,...,n_1$;
        \item $\lambda_i(G_2)+m$, for $i=2,3,...,n_2$;
        \item $\frac{(n_1+k+n_2+2)\pm \sqrt{(n_1+k-n_2-2)^2+4(2k-\lambda_i(G))}}{2}$ , for $i=2,3,...,n$;
        \item $n_2+2$, repeated $m-n$ times;
        \item all the roots of the following equation\\
        $\lambda(\lambda^3-(m+n+n_1+k+n_2+2)\lambda^2+(2m+(n_1+k+n)(m+n_2+2)+nk-2k)\lambda+2(n+m)k-2(n_1+k+n)m-nk(m+n_2+2))=0$.
    \end{enumerate}

  \begin{proof}
With a suitable labeling of the vertices of $G^S\vee{(G^\bullet_1,
G^\circ_2)}$, we can write the Laplacian matrix of
$G^S\vee{(G^\bullet_1, G^\circ_2)}$ as
$$ L(G^S\vee{(G^\bullet_1 , G^\circ_2)})=
    \begin{pmatrix}
    (n_1+k)I_n & -M& -J_{n\times n_1}& \textbf{0}_{n\times n_2}\\
    -M^T& (n_2+2)I_m & \textbf{0}_{m\times n_1}& -J_{m\times n_2}\\
    -J_{n_1\times n}& \textbf{0}_{n_1\times m} & L(G_1)+nI_{n_1} & \textbf{0}_{n_1\times n_2}\\
    \textbf{0}_{n_1\times n}& -J_{n_2\times m}& \textbf{0}_{n_2\times n_1}& L(G_2)+mI_{n_2}\\
    \end{pmatrix},
    $$
where $M$ denotes vertex-edge incidence matrix of $G$ and $I_n$ the
identity matrix of order $n$. By comparing the Laplacian matrix
$L(G^S\vee{(G^\bullet_1 , G^\circ_2)})$ with the double join matrix
$\mathscr{D}_j$, we take $p=n$, $q=m$, $r=n_1$, $s=n_2$, $c=-1$ and
$A=(n_1+k)I_n$, $B=-M$, $C=(n_2+2)I_m$, $D=L(G_1)+nI_{n_1}$,
$E=L(G_2)+mI_{n_2}$ in Theorem 3. Since
$MM^T=\left|L\right|\left(G\right)=kI_n+A(G)=2kI_n-L(G)$,
$L(G_1)\textbf{1}_{n_1}=0\textbf{1}_{n_1}$ and
$L(G_2)\textbf{1}_{n_2}=0\textbf{1}_{n_2}$. Then we have
    \begin{itemize}
        \item $a_i=n_1+k$ for $i=1,2,...,n$,\;\; $a_1=n_1+k$;
        \item $b_i^2=2k-\lambda_i(G)$ for $i=1,2,...,n$,\;\; $b_1^2=2k$;
        \item $c_i=n_2+2$ for $i=1,2,...,m$;
        \item $d_i=\lambda_i(G_1)+n$ for $i=1,2,...,n_1$, \;\;$d_1=n$;
        \item $e_i=\lambda_i(G_2)+m$ for $i=1,2,...,n_2$,\;\; $e_1=m$.
    \end{itemize}
Now plugging these values into Theorem 3, we obtain the required
result.
  \end{proof}
\paragraph{Remark 5.} Remark that the subdivision double join $G^S\vee{(G^\bullet_1
, G^\circ_2)}$ becomes the subdivision-vertex join defined in
\cite{Indulal2012} whenever $G_2$ is a null graph. Similarly, the
subdivision double join $G^S\vee{(G^\bullet_1 , G^\circ_2)}$ becomes
the subdivision-edge join defined in \cite{Indulal2012} whenever
$G_1$ is a null graph. In \cite{Liu2017}, Liu and Zhang determined
the Laplacian spectra of subdivision-vertex join and
subdivision-edge join. Clearly, Theorem 4 generalizes the results of
both Theorems 2.7 and 3.4 in \cite{Liu2017}.\\

Next, we give a complete description of the Laplacian spectra of the
$Q$-graph double join $G^Q\vee{(G^\bullet_1, G^\circ_2)}$ for a
regular graph $G$ and two arbitrary graphs $G_1$, $G_2$.

\paragraph{Theorem 6.} Let $G$ be a $k$-regular graph with $n$ vertices and $m$
edges. Also let $G_1$ and $G_2$ be two arbitrary graph with $n_1$
and $n_2$ vertices, respectively. Then the Laplacian spectrum of
$G^Q\vee{(G^\bullet_1, G^\circ_2)}$ consists of:
\begin{enumerate}[(i)]
  \item $\lambda_i(G_1)+n$, for $i=2,3,...,n_1$;
  \item $\lambda_i(G_2)+m$, for $i=2,3,...,n_2$;
  \item $\frac{(n_1+k+n_2+2+\lambda_i(G))\pm \sqrt{(n_1+k-n_2-2-\lambda_i(G))^2+4(2k-\lambda_i(G))}}{2}$ ,  for $i=2,3,...,n$;
  \item $n_2+2k+2$, repeated $m-n$ times;
  \item four roots of the equation\\
  $\lambda(\lambda^3-(m+n+n_1+k+n_2+2)\lambda^2+(2m+(n_1+k+n)(m+n_2+2)+nk-2k)\lambda+2(n+m)k-2(n_1+k+n)m-nk(m+n_2+2))=0.$
\end{enumerate}
\begin{proof}

With a proper labeling of vertices, the Laplacian matrix of
$G^Q\vee{(G^\bullet_1, G^\circ_2)}$ can be written as
    $$L(G^Q\vee{(G^\bullet_1 , G^\circ_2)})=
    \begin{pmatrix}
    (n_1+k)I_n & -M& -J_{n\times n_1}& \textbf{0}_{n\times n_2}\\
    -M^T& (n_2+2)I_m+L(l(G)) & \textbf{0}_{m\times n_1}& -J_{m\times n_2}\\
    -J_{n_1\times n}& \textbf{0}_{n_1\times m} & L(G_1)+nI_{n_1} & \textbf{0}_{n_1\times n_2}\\
    \textbf{0}_{n_1\times n}& -J_{n_2\times m}& \textbf{0}_{n_2\times n_1}& L(G_2)+mI_{n_2}\\
    \end{pmatrix},$$
where $l(G)$ denotes the line graph of $G$.

Now, comparing the Laplacian matrix $L(G^Q\vee{(G^\bullet_1 ,
G^\circ_2)})$ with the double join matrix $\mathscr{D}_j$, we take
$A=(n_1+k)I_n$, $B=-M$, $C=(n_2+2)I_m+L(l(G))$, $D=L(G_1)+nI_{n_1}$,
$E=L(G_2)+mI_{n_2}$ in Theorem 3. Since $D(l(G))=2(k-1)I_m$,
$A(l(G))=M^TM-2I_m$. Then $L(l(G))=2kI_m-M^TM$, which implies that
$C=(n_2+2k+2)I_m-M^TM$. Note that $MM^T=2kI_n-L(G)$. Since $MM^T$
and $M^TM$ have same nonzero eigenvalues. Then the spectrum of $C$
consists of: $c_i=n_2+\lambda_i(G)+2$ for $i=1,2,...,n$ and
$c_j=n_2+2k+2$ for $j=n+1,n+2,...,m$. Furthermore, it follows from
Theorem 3.38 in \cite{Cvetkovic1980} that
$C\textbf{Y}_j=c_j\textbf{Y}_j$ satisfy $M\textbf{Y}_j=\textbf{0}_n$
for $j=n+1,n+2,...,m$. Therefore, one has
\begin{itemize}
    \item $a_i=n_1+k$ for $i=1,2,...,n$,\;\; $a_1=n_1+k$;
    \item $b_i^2=2k-\lambda_i(G)$ for $i=1,2,...,n$,\;\; $b_1^2=2k$;
    \item $c_i=n_2+\lambda_i(G)+2$ for $i=1,2,...,n$,\;\;$c_j=n_2+2k+2$ for $j=n+1,...,m$;
    \item $d_i=\lambda_i(G_1)+n$ for $i=1,2,...,n_1$, \;\;$d_1=n$;
    \item $e_i=\lambda_i(G_2)+m$ for $i=1,2,...,n_2$, \;\;$e_1=m$.
\end{itemize}
Now from Theorem 3, the required result follows.
\end{proof}

The following result describes the Laplacian spectra of the
$R$-graph double join $G^R\vee{(G^\bullet_1, G^\circ_2)}$ for a
regular graph $G$ and two arbitrary graphs $G_1$, $G_2$.

\paragraph{Theorem 7.} Let $G$ be a $k$-regular graph with $n$ vertices and $m$
edges. Let $G_1$ and $G_2$ be two arbitrary graph witn $n_1$ and
$n_2$ vertices, respectively. Then the Laplacian spectrum of
$G^R\vee{(G^\bullet_1, G^\circ_2)}$ consists of:
     \begin{enumerate}[(i)]
        \item $\lambda_i(G_1)+n$, for $i=2,3,...,n_1$;
        \item $\lambda_i(G_2)+m$, for $i=2,3,...,n_2$;
        \item $\frac{(n_1+\lambda_i(G)+k+n_2+2)\pm \sqrt{(n_1+\lambda_i(G)+k-n_2-2)^2+4(2k-\lambda_i(G))}}{2}$ ,  for $i=2,3,...,n$;
        \item $n_2+2$, repeated $m-n$ times;
        \item four roots of the equation\\
        $\lambda(\lambda^3-(m+n+n_1+k+n_2+2)\lambda^2+(2m+(n_1+k+n)(m+n_2+2)+nk-2k)\lambda+2(n+m)k-2(n_1+k+n)m-nk(m+n_2+2))=0.$
     \end{enumerate}
\begin{proof}
With a proper labeling of vertices, the Laplacian matrix of
$G^R\vee{(G^\bullet_1, G^\circ_2)}$ can be written as
$$L(G^R\vee{(G^\bullet_1 , G^\circ_2)})=
\begin{pmatrix}
L(G)+(n_1+k)I_n & -M& -J_{n\times n_1}& \textbf{0}_{n\times n_2}\\
-M^T& (n_2+2)I_m & \textbf{0}_{m\times n_1}& -J_{m\times n_2}\\
-J_{n_1\times n}& \textbf{0}_{n_1\times m} & L(G_1)+nI_{n_1} & \textbf{0}_{n_1\times n_2}\\
\textbf{0}_{n_1\times n}& -J_{n_2\times m}& \textbf{0}_{n_2\times n_1}& L(G_2)+mI_{n_2}\\
\end{pmatrix}.
$$
Using the same technique as the proof of Theorem 4, we obtain
\begin{itemize}
  \item $a_i=n_1+k+\lambda_i(G)$ for $i=1,2,...,n$,\;\; $a_1=n_1+k$;
  \item $b_i^2=2k-\lambda_i(G)$ for $i=1,2,...,n$,\;\; $b_1^2=2k$;
  \item $c_i=n_2+2$ for $i=1,2,...,m$;
  \item $d_i=\lambda_i(G_1)+n$ for $i=1,2,...,n_1$,\;\; $d_1=n$;
  \item $e_i=\lambda_i(G_2)+m$ for $i=1,2,...,n_2$,\;\; $e_1=m$.
\end{itemize}
Now plugging these values in Theorem 3, we obtain the desired
result.
\end{proof}

For the total double join $G^T\vee{(G^\bullet_1, G^\circ_2)}$, we
describe the Laplacian spectra in the following results.

\paragraph{Theorem 8.}Let $G$ be a $k$-regular graph with $n$ vertices and $m$
edges. Let $G_1$ and $G_2$ be two arbitrary graph with $n_1$ and
$n_2$ vertices, respectively. Then the Laplacian spectrum of
$G^T\vee{(G^\bullet_1, G^\circ_2)}$ consists of:
     \begin{enumerate}[(i)]
        \item $\lambda_i(G_1)+n$, for $i=2,3,...,n_1$;
        \item $\lambda_i(G_2)+m$, for $i=2,3,...,n_2$;
        \item $\frac{(n_1+2\lambda_i(G)+k+n_2+2)\pm \sqrt{(n_1+k-n_2-2)^2+4(2k-\lambda_i(G))}}{2}$ ,  for $i=2,3,...,n$;
        \item $n_2+2k+2$, repeated $m-n$ times;
        \item four roots of the equation\\
        $\lambda(\lambda^3-(m+n+n_1+k+n_2+2)\lambda^2+(2m+(n_1+k+n)(m+n_2+2)+nk-2k)\lambda+2(n+m)k-2(n_1+k+n)m-nk(m+n_2+2))=0$.
     \end{enumerate}
\begin{proof}
The Laplacian matrix of $G^T\vee{(G^\bullet_1, G^\circ_2)}$ can be
expressed as follows
$$
\begin{pmatrix}
L(G)+(n_1+k)I_n & -M& -J_{n\times n_1}& \textbf{0}_{n\times n_2}\\
-M^T& L(l(G))+(n_2+2)I_m & \textbf{0}_{m\times n_1}& -J_{m\times n_2}\\
-J_{n_1\times n}& \textbf{0}_{n_1\times m} & L(G_1)+nI_{n_1} & \textbf{0}_{n_1\times n_2}\\
\textbf{0}_{n_1\times n}& -J_{n_2\times m}& \textbf{0}_{n_2\times n_1}& L(G_2)+mI_{n_2}\\
\end{pmatrix}.
$$
Using the similar technique to the proof of Theorem 6, we obtain
         \begin{itemize}
            \item $a_i=n_1+k+\lambda_i(G)$ for $i=1,2,...,n$,\;\; $a_1=n_1+k$;
            \item $b_i^2=2k-\lambda_i(G)$ for $i=1,2,...,n$, \;\;$b_1^2=2k$;
            \item $c_i=n_2+\lambda_i(G)+2$ for $i=1,2,...,n$,\;\;$c_j=n_2+2k+2$ for $j=n+1,...,m$;
            \item $d_i=\lambda_i(G_1)+n$ for $i=1,2,...,n_1$,\;\; $d_1=n$;
            \item $e_i=\lambda_i(G_2)+m$ for $i=1,2,...,n_2$, \;\;$e_1=m$.
         \end{itemize}
Now substituting these values in Theorem 3, we get the expected
result.
\end{proof}

\paragraph{Remark 9.} If $G_2$ is a null graph, then our $Q$-graph double
join($R$-graph double join, total double join) reduces to $Q$-graph
vertex join($R$-graph vertex join\cite{Liu2015}, total vertex join,
respectively). Similarly, If $G_1$ is a null graph, then our
$Q$-graph double join($R$-graph double join, total double join)
reduces to $Q$-graph edge join($R$-graph edge join\cite{Liu2015},
total edge join, respectively). Then Theorems 6, 7 and 8 can help us
to determine completely Laplacian spectra of these join operations
of graphs.

\section*{4. Conclusion}

\indent \indent Here we introduce the conception of double join
matrix and provide a complete description about its eigenvalues and
the corresponding eigenvectors. Further, applying the result
obtained for the double join matrix, we give an explicit complete
characterization of the Laplacian spectra of four variants of double
join operations of graphs in terms of the Laplacian spectra of the
factor graphs. These results not only generalize some well-known
results, but also describe completely the eigenvectors corresponding
to all the Laplacian eigenvalues of these graphs.

As described in the Introduction, many families of pairs of
cospectral graphs may be constructed by using some graph operations.
Assume that $G$ and $H$ (not necessarily distinct) are two Laplacian
cospectral regular graphs, $G_1$ is Laplacian cospectral with
$H_1$(not necessarily distinct) and $G_2$ is Laplacian cospectral
with $H_2$(not necessarily distinct). Then
$G^S\vee{(G^\bullet_1,G^\circ_2)}$ and
$H^S\vee{(H^\bullet_1,H^\circ_2)}$ are Laplacian cospectral.
Similarly, we can also construct many families of pairs of Laplacian
cospectral graphs for other variants of double join operations.

The degree Kirchhoff index and the number of spanning trees of some
graph operations have been studied extensively(see Introduction).
Our results can also help us to compute the number of spanning trees
and Kirchhoff index for four variants of double join operations of
graphs.

Before the end of this paper, we see easily that the Laplacian
matrix of join graph is also a double join matrix by choosing
$B=\textbf{0}_{p\times q}$, $C=\textbf{0}_{q\times q}$ and $s=0$ in
the double join matrix $\mathscr{D}_j$. Thus the nonzero Laplacian
eigenvalues of the join graph can be obtained from the corresponding
eigenvalues of double join matrix. Hence our result also generalizes
the classical result about the Laplacian spectrum of the usual join
graph obtained in \cite{Merris1998}.
 \\
\\
\textbf{Acknowledgements} This work was in part supported by NNSFC
(Nos. 11371328, 11671053) and by the Natural Science Foundation of
Zhejiang Province, China (No. LY15A010011).

\end{document}